\documentclass[12pt]{amsart}
\usepackage[a4paper,bottom=2cm,top=1cm,left=2cm,right=1.3cm]{geometry}
\usepackage[english]{babel}
\usepackage{amsfonts,amsmath,latexsym,amssymb,amsthm,mathrsfs,upref,mathptmx} 

\newtheorem{theorem}{Theorem}

\theoremstyle{definition}

\begin{document}

\title{Multi-interval Sturm--Liouville boundary-value problems \\with distributional potentials}

\author{Andrii Goriunov}

\address{Institute of Mathematics of National Academy of Sciences of Ukraine, Kyiv, Ukraine}
\email{goriunov@imath.kiev.ua}



%\thanks{This research is supported by the grant no. 03-01-12 of the National Academy of Sciences of Ukraine 
%(under the joint Ukrainian-Russian project of the NAS of Ukraine and the Siberian Branch of RAS).}

\begin{abstract}
We study the multi-interval boundary-value Sturm-Liouville problems with distributional potentials. 
For the corresponding symmetric operators boundary triplets are found and the constructive descriptions 
of all self-adjoint, maximal dissipative and maximal accumulative extensions and generalized resolvents 
in terms of homogeneous boundary conditions are given. 
It is shown that all real maximal dissipative and maximal accumulative extensions are self-adjoint 
and all such extensions are described.
\end{abstract}

\maketitle

In recent years the interest in multi-interval differential and quasi-differential operators has increased 
(see \cite{EvZ_86,EvZ,Sok_03,Sok_06}). 
The main attention is paid to the case where a (quasi-)differential expression is formally self-adjoint.
From the operator-theoretic point of view this corresponds to the situation where 
we investigate extensions of a symmetric (quasi-)differential operator with equal deficiency indices 
in the direct sum of Hilbert spaces on the basis of Glazman-Krein-Naimark theory 
\cite{Z_75,EM-book,Z_book,Naim_book}. 
In the present paper we develop another approach to such problems 
based on the concept of boundary triplets \cite{Gorb_book,Koch_mzm_79}.

Let $m\in\mathbb{N}$, $a=a_0<a_1<\dots<a_m=b$ be a partition of a finite interval $[a,b]$ into $m$ parts  
and on every interval $(a_{i-1},a_{i})$, $i \in \{1,\dots, m\}$ 
let the formal Sturm-Liouville expression 
\begin{equation}\label{S-L_expr_i}
l_i(y) = -(p_i(t)y')'+q_i(t)y
\end{equation}
be given.
Here, the measurable finite functions $p_i$ and $Q_i$ are such that  
\begin{equation}\label{GM_cond_i}
1/p_i, \,\, Q_i/p_i,\,\,  Q_i^2/p_i \in L_1\left([a_{i-1},a_{i}],\mathbb{R}\right),
\end{equation}
the potentials $q_i = Q_i'$ and the derivative is understood in the sense of distributions. 

For $m=1$ the boundary-value problems for the formal differential expression \eqref{S-L_expr_i} 
under assumptions \eqref{GM_cond_i} 
were investigated in \cite{GM_MFAT_10} on the basis of its regularization by Shin--Zettl quasiderivatives. 
In this paper the most of the results of \cite{GM_MFAT_10} is extended onto the case of an arbitrary $m\in\mathbb{N}$.

We introduce the quasi-derivatives 
\begin{align*}
&D_i^{[0]} y = y, \\
&D_i^{[1]} y = p_iy' - Q_iy, \\
&D_i^{[2]} y = (D_i^{[1]} y)' + {Q_i\over p_i}D_i^{[1]} y + {Q_i^2\over p_i}y.
\end{align*}
on every interval $(a_{i-1},a_{i})$, as in \cite{GM_MFAT_10}.

Then the maximal and minimal operators 
\begin{align*}
&L_{i,1}:y \to
l_i[y],\quad \text{Dom}(L_{i,1}) := 
\left\{y \in L_2 \left| y,\, D_i^{[1]}y \in AC([a_{i-1},a_i],\mathbb{C}), \,\,D_i^{[2]} y \in L_2\right.\right\},\\
&L_{i,0}:y \to
l_i[y],\quad \text{Dom}(L_{i,0}) := 
\left\{y \in \text{Dom}(L_{i,1}) \left| D_i^{[k]}y(a_{i-1}) = D_i^{[k]}y(a_i) = 0, 
\,\,k = 0,1\right.\right\}
\end{align*}
are defined in the spaces $L_2\left((a_{i-1},a_{i}),\mathbb{C}\right)$.
According to \cite{GM_MFAT_10} the operators $L_{i,1}$, $L_{i,0}$ are closed and densely defined in  
$L_2\left([a_{i-1},a_{i}],\mathbb{C}\right)$.
The operator $L_{i,0}$ is symmetric with the deficiency indices $\left({2,2} \right)$ and 
$$L_{i,0}^* = L_{i,1},\quad L_{i,1}^* = L_{i,0}.$$

Recall that a \emph{boundary triplet} of a closed densely defined symmetric operator $T$ 
with equal (finite or infinite) deficiency indices is called a triplet 
$\left(H, \Gamma_1 ,\Gamma_2 \right)$ where $H$ is an auxiliary Hilbert space and
$\Gamma_1$, $\Gamma_2$ are the linear maps from $\text{Dom}(T^*)$ to $H$ such that  
\begin{enumerate}
\item for any $f,g \in \text{Dom}\left(T^*\right)$ there holds
\[
\left(T^*f,g\right)_\mathcal{H} - \left(f,T^*g\right)_\mathcal{H} = \left(\Gamma_1f,\Gamma_2g\right)_H  -
\left(\Gamma_2f,\Gamma_1g\right)_H;
\]
\item for any $g_1, g_2 \in H$ there is a vector $f\in \text{Dom}\left(T^*\right)$ 
such that $\Gamma_1 f = g_1$ and $\Gamma_2 f = g_2$.
\end{enumerate}

It is proved in \cite{GM_MFAT_10} that for every $i=1,\dots, m$ 
the triplet $(\mathbb{C}^{2}, \Gamma_{1,i}, \Gamma_{2,i})$, 
where $\Gamma_{1,i}, \Gamma_{2,i}$ are linear maps
%\begin{equation} \label{PGZ_i}
\[ \Gamma_{1,i}y := \left( D_i^{[1]}y(a_{i-1}+), -D_i^{[1]}y(a_{i}-)\right), \, 
	\Gamma_{2,i}y := \left( y(a_{i-1}+), y(a_{i}-)\right), \]
	%\end{equation}
from $\text{Dom}(L_{i,1})$ to $\mathbb{C}^{2}$ 
is a boundary triplet for the operator $L_{i,0}$.

We consider the space $L_2\left([a,b],\mathbb{C}\right)$ as a direct sum 
$\oplus_{i=1}^mL_2\left([a_{i-1},a_{i}],\mathbb{C}\right)$ which consists of vector functions
$f= \oplus_{i=1}^m f_i$ such that $f_i \in L_2\left([a_{i-1},a_{i}],\mathbb{C}\right)$.
In this space we consider operators $L_{\text{max}}=\oplus_{i=1}^mL_{i,1}$ and 
$L_{\text{min}}=\oplus_{i=1}^mL_{i,0}$.

%\begin{lemma}\label{lem_maxmin}
Then the operators $L_{\text{max}}$, $L_{\text{min}}$ are closed and densely defined in 
$L_2\left([a,b],\mathbb{C}\right)$.
The operator $L_{\text{min}}$ is symmetric with the deficiency indices $\left({2m,2m} \right)$ and 
$$L_{\text{min}}^* = L_{\text{max}},\quad L_{\text{max}}^* = L_{\text{min}}.$$
%\end{lemma}

Note that the minimal operator $L_{\text{min}}$ may be not semi-bounded  even 
in the case of a single-interval boundary-value problem 
since the function $p$ may reverse sign.

\begin{theorem}\label{th_PGZ}
The triplet $(\mathbb{C}^{2m}, \Gamma_{1}, \Gamma_{2})$ 
where $\Gamma_{1}, \Gamma_{2}$ are linear maps 
%\begin{equation} \label{PGZ}
\[ \Gamma_{1}y := \left( \Gamma_{1,1}y, \Gamma_{1,2}y,\dots, \Gamma_{1,m}y\right), \, 
	\Gamma_{2}y := \left( \Gamma_{2,1}y, \Gamma_{2,2}y,\dots, \Gamma_{2,m}y\right) \]
%\end{equation}
from $\text{Dom}(L_{\text{max}})$ onto $\mathbb{C}^{2m}$ is a boundary triplet for $L_{\text{min}}$.
\end{theorem}

Denote by $L_K$ the restriction of $L_{\text{max}}$ onto the set of functions 
${y(t) \in \text{Dom}(L_{\text{max}})}$ satisfying the homogeneous boundary condition 
%\begin{equation} \label{L_K}
 \[\left( {K - I} \right)\Gamma_1 y + i\left( {K + I} \right)\Gamma_2 y = 0.\]
%\end{equation}
%where $K$ is a bounded operator in $\mathbb{C}^{2m}$.

Similarly, denote by $L^K$ the restriction of $L_{\text{max}}$ onto the set of functions
${y(t) \in \text{Dom}(L_{\text{max}})}$ satisfying the homogeneous boundary condition 
%\begin{equation} \label{L^K}
 \[\left( {K - I} \right)\Gamma_1 y - i\left( {K + I} \right)\Gamma_2 y = 0.\]
%\end{equation}
Here $K$ is a bounded operator in $\mathbb{C}^{2m}$.

The constructive description of the various classes of extensions of the operator $L_{\text{min}}$ is given 
by the following theorem.

\begin{theorem}\label{th_ext}
Every $L_K$ with $K$ being a contracting operator in $\mathbb{C}^{2m}$ is a maximal dissipative extension of $L_{\text{min}}$.
Similarly every $L^K$ with $K$ being a contracting operator in $\mathbb{C}^{2m}$ 
is a maximal accumulative extension of the operator $L_{\text{min}}$.

Conversely, for any maximal dissipative (respectively, maximal accumulative) extension $\widetilde{L}$ 
of the operator $L_{\text{min}}$ there exists the unique contracting operator $K$ such that 
$\widetilde{L} = L_K$\,\, (respectively, $\widetilde{L} = L^K$). 

The extensions $L_K$ and $L^K$ are self-adjoint if and only if $K$ is a unitary operator 
on $\mathbb{C}^{2m}$. 
\end{theorem}

\vspace*{2mm}

%\begin{definition}
Recall that a linear operator $T$ acting in $L_2([a,b], \mathbb{C})$ is called \emph{real} if:
\begin{enumerate}
\item   for every function $f$ from $\text{Dom}(T)$ the complex conjugate function $\overline{f}$ 
also lies in $\text{Dom}(T)$;
\item the operator $T$ maps complex conjugate functions into complex conjugate functions, that is 
$T(\overline{f}) = \overline{T(f)}$.
\end{enumerate}
%\end{definition}

One can see that the maximal and minimal operators are real.

\begin{theorem}\label{th_real}
All real maximal dissipative and maximal accumulative extensions of the minimal operator $L_{\text{min}}$ 
are self-adjoint.
The self-adjoint extension $L_K$ or $L^K$ is real if and only if the unitary matrix $K$ is symmetric.
\end{theorem}

\vspace{2mm}

%\begin{definition}
Let us recall that a \emph{generalized resolvent} of a closed symmetric operator $T$ 
in a Hilbert space $\mathcal{H}$ is an operator-valued function $\lambda\mapsto R_\lambda$ 
defined on $\mathbb{C} \setminus \mathbb{R}$, which can be represented as
\[
R_\lambda f = P^+ \left( T^+ - \lambda I^+\right)^{- 1}f, \quad f \in \mathcal{H},
\]
where $T^+$ is a self-adjoint extension of $T$ 
which acts in a certain Hilbert space $\mathcal{H}^+\supset\mathcal{H}$,
$I^+$ is the identity operator on $\mathcal{H}^+$, 
and $P^+$ is the orthogonal projection operator from $\mathcal{H}^+$ onto $\mathcal{H}$.
It is known that an operator-valued function $R_\lambda$ % \quad (\text{Im}\lambda \neq 0)$
is a generalized resolvent of a symmetric operator $T$ if and only if it can be represented as
\[
\left( R_\lambda f, g \right)_\mathcal{H} = \int_{-\infty}^{+\infty}\frac{d\left(F_\mu f, g\right)}{\mu - \lambda},
  \quad f, g \in \mathcal{H},
\] 
where $F_\mu$ is a generalized spectral function of the operator $T$, i. e.  
$\mu\mapsto F_\mu$ is an operator-valued function $F_\mu$ defined on $\mathbb{R}$ and taking
values in the space of continuous linear operators in $\mathcal{H}$ with the following properties:
\begin{enumerate}
\item  for $\mu_2 > \mu_1$ the difference $F_{\mu_2} - F_{\mu_1}$
 is a bounded non-negative operator;
\item $F_{\mu +} = F_\mu$ for any real $\mu$;
\item for any $x \in \mathcal{H}$ there holds
\[ \lim\limits_{\mu \rightarrow - \infty}^{}||F_\mu x ||_\mathcal{H} = 0,
 \quad \lim\limits_{\mu \rightarrow + \infty}^{} ||{F_\mu x - x} ||_\mathcal{H} = 0.\]
\end{enumerate}

The following theorem provides a description of all generalized resolvents of the operator $L_{\text{min}}$.
%\end{definition}

\begin{theorem}\label{th_gen_res}
$1)$ %Нехай $\lambda$ -- комплексне число, $\text{Im}\lambda < 0$.
Every generalized resolvent $R_\lambda$ of the operator $L_{\text{min}}$
in the half-plane $\text{Im}\lambda < 0$ acts by the rule $R_\lambda h = y$, 
where $y$ is a solution of the boundary-value problem 
$$ l(y) = \lambda y + h,$$
$$
 \left( {K(\lambda) - I} \right)\Gamma _{1} f + i\left( {K(\lambda) + I} \right)\Gamma _{2} f = 0.
$$
Here $h(x) \in L_2([a,b], \mathbb{C})$ and
$K(\lambda)$ is a $2m\times 2m$ matrix-valued function which is holomorphic 
in the lower half-plane and satisfies $||K(\lambda)|| \leq 1$.

$2)$ In the half-plane $\text{Im}\lambda > 0$ every generalized resolvent of $L_{\text{min}}$ 
acts by the rule $R_\lambda h = y$ where
$y$ is a solution of the boundary-value problem 
%де
%Існує взаємно однозначна відповідність між узагальненими
%резольвентами оператора $L_{\text{min}}$ і крайовими задачами
$$ l(y) = \lambda y + h,$$
$$
 \left( {K(\lambda) - I} \right)\Gamma _{1} f - i\left( {K(\lambda) + I} \right)\Gamma _{2} f = 0.
$$
Here $h(x) \in L_2([a,b], \mathbb{C})$ and $K(\lambda)$ is a $2m\times 2m$ matrix-valued function 
which is holomorphic in the lower half-plane and satisfies $||K(\lambda)|| \leq 1$.
\vskip 2mm
The parametrization of the generalized resolvents by the matrix-valued functions $K$ is bijective.
\end{theorem}

\medskip

%In the present paper the operators corresponding to multi-interval boundary-value Sturm-Liouville problems 
%with distributional potentials are studied. 
%For these symmetric operators boundary triplets are found and constructive descriptions 
%of all self-adjoint, maximal dissipative and maximal accumulative extensions and generalized resolvents 
%in terms of homogeneous boundary conditions are given. 
%It is shown that all real maximal dissipative and maximal accumulative extensions are self-adjoint 
%and all such extensions are described.
%
%\medskip

{

\end{document}